\documentclass[12pt]{article}
\usepackage{amssymb}

\newtheorem{theorem}{Theorem}
\newtheorem{proposition}[theorem]{Proposition}
\newtheorem{corollary}[theorem]{Corollary}

\newenvironment{proof}[1][Proof]{\textbf{#1.} }{\ \rule{0.5em}{0.5em}}

\begin{document}

\title{Some basic properties of Lagrangians}
\author{\textsc{Ioan Bucataru} \\
{\small Faculty of Mathematics, ``Al.I.Cuza'' University,} \\
{\small Ia\c{s}i, 700506, Romania, email: bucataru@uaic.ro}}
\date{}
\maketitle

\noindent\textbf{Abstract} {\small Consider $L$ a regular
Lagrangian, $S$ the canonical semispray, and $h$ the horizontal
projector of the canonical nonlinear connection. We prove that if
the Lagrangian is constant along the integral curves of the
Euler-Lagrange equations then it is constant along the horizontal
curves of the canonical nonlinear connection. In other words
$S(L)=0$ implies $d_hL=0$. If the Lagrangian $L$ is homogeneous of
order $k\neq 1$ then $L$ is a conservation law and hence $d_hL=0$.
We give an example of nonhomogeneous Lagrangians for which
$d_hL\neq 0$.}

\vspace*{3mm}



\section*{Introduction}

For a Lagrangian $L$ its energy $E_L=\mathbb{C}(L)-L$ is constant
along the solution curves of the Euler-Lagrange equations. In
other words, we can say that the energy $E_L$ is constant along
the integral curves of the canonical semispray $S$, and this
property can be written as $S(E_L)=0$. For the particular case
when the Lagrangian $L$ is homogeneous of order $k\neq 1$, with
respect to the fibre coordinates, its energy is proportional with
the Lagrangian, which means that $E_L=(k-1)L$. If this is the
case, the horizontal curves of the canonical nonlinear connection
coincide with the integral curves of the canonical semispray.
Hence, the Lagrangian is constant along the horizontal curves of
the canonical nonlinear connection. We can express this property
as $d_hL=0$, where $h$ is the horizontal projector that
corresponds to the canonical nonlinear connection. a similar
result has been obtained by J. Grifone in \cite{[grifone]} and J.
Szenthe in \cite{[szenthe2]}.

In two papers \cite{[szenthe1]} and \cite{[szilasi]} it is stated
that for any Lagrangian L we have that $d_hL=0$, which is not true
if $L$ is not homogeneous of order one. For nonhomogeneous
Lagrangians the integral curves of the canonical semispray do not
coincide with the horizontal curves of the nonlinear connection.
We prove that if the Lagrangian is constant along the integral
curves of the canonical semispray then it is constant also along
the horizontal curves of the canonical nonlinear connection. In
other words we have that $S(L)=0$ implies $d_h(L)=0$. In general
for nonhomogeneous Lagrangian we have that $S(L)\neq 0$ and
$d_hL\neq 0$ and we give examples of such Lagrangians. However, we
don't know if the converse implication, $d_hL=0$ implies $S(L)=0$,
is true.

\section{Geometric structures on tangent bundle}

In this section we introduce the geometric objects we are going to
deal with in this paper such as: Liouville vector field, tangent
structure, semispray and nonlinear connection. Since most of them
live on the total space of the tangent bundle of a differentiable
manifold, we briefly recall some basic properties of the tangent
bundle.

Let $M$ be a real, $n$-dimensional manifold of $C^{\infty}$-class
and denote by $(TM,\pi,M)$ its tangent bundle. We denote by
$\widetilde{TM}=TM\setminus 0$ the tangent bundle with zero
section removed. If $(U, \phi=(x^i))$ is a local chart at $p\in M$
from a fixed atlas of $C^{\infty}$-class, then we denote by
$(\pi^{-1}(U), \Phi=(x^i, y^i))$ the induced local chart at $u\in
\pi^{-1}(p) \subset TM$. The linear map $\pi_{*,u}:T_uTM
\rightarrow T_{\pi(u)}M$, induced by the canonical submersion
$\pi$, is an epimorphism of linear spaces for each $u\in TM$.
Therefore, its kernel determines a regular, $n$-dimensional,
integrable distribution $V:u\in TM \mapsto V_uTM:=\mathrm{Ker}
\pi_{*,u} \subset T_uTM$, which is called the \textit{vertical
distribution}. For every $u \in TM$, $\{{\partial}/{\partial
y^i}|_u\}$ is a basis of $V_uTM$, where $\{{\partial}/{\partial
x^i}|_u, {\partial}/{\partial y^i}|_u\}$ is the natural basis of
$T_uTM$ induced by a local chart. Denote by $\mathcal{F}(TM)$ the
ring of real-valued functions over $TM$ and by $\mathcal{X}(TM)$
the $\mathcal{F}(TM)$-module of vector fields on $TM$. We also
consider $\mathcal{X}^v(TM)$ the $\mathcal{F}(TM)$-module of
vertical vector fields on $TM$. An important vertical vector field
is $\mathbb{C} = y^i({\partial}/{\partial y^i})$, which is called
the \textit{Liouville vector field}.

The mapping $J: \mathcal{X}(TM) \rightarrow \mathcal{X}(TM)$ given
by $J=({\partial}/{\partial y^i}) \otimes dx^i $ is called the
\textit{tangent structure} and it has the following properties:
Ker $J$ = Im $J$ = $\mathcal{X}^v(TM)$; rank $J=n$ and $J^2=0$.
The cotangent structure is defined as $J^*=dx^i \otimes
({\partial}/{\partial y^i})$ and has similar properties.

A vector field $S\in\chi(TM)$, which is differentiable of
$C^{\infty}$-class on $\widetilde{TM}$ and only continuous on the
null section, is called a \textit{semispray}, or a second order
vector field, if $JS=\mathbb{C}$. In local coordinates a semispray
can be represented as follows:
\begin{equation}
S=y^i\frac{\partial}{\partial x^i}-2G^i(x,y)
\frac{\partial}{\partial y^i}. \label{smispray}
\end{equation} We refer to the functions $G^i(x,y)$ as to the local coefficients
of the semispray $S$.

A \textit{nonlinear connection} $N$ on $TM$ is an $n$-dimensional
distribution, which is also called the horizontal distribution,
$N: u\in TM \mapsto N_uTM\subset T_uTM$ that is supplementary to
the vertical distribution $VTM$. This means that for every $u\in
TM$ we have the direct decomposition:
\begin{equation}
T_uTM=N_uTM \oplus V_uTM. \label{thv} \end{equation} We denote by
$h$ and $v$ the horizontal and the vertical projectors that
correspond to the above decomposition and by $\mathcal{X}^h(TM)$
the $\mathcal{F}(TM)$-module of horizontal vector fields on $TM$.
For every $u=(x,y)\in TM$ we denote by ${\delta}/{\delta
x^i}|_u=h({\partial}/{\partial x^i}|_u).$ Then $\{{\delta}/{\delta
x^i}|_u, {\partial}/{\partial y^i}|_u\}$ is a basis of $T_uTM$
adapted to the decomposition (\ref{thv}). With respect to the
natural basis $\{{\partial}/{\partial x^i}|_u,
{\partial}/{\partial y^i}|_u\}$ of $T_uTM$, the horizontal
components of the adapted basis have the expression:
\begin{equation}
\left.\frac{\delta}{\delta x^i}\right|_u =
\left.\frac{\partial}{\partial x^i}\right|_u -
N^j_i(u)\left.\frac{\partial}{\partial y^j}\right|_u, \ u \in TM.
\label{coefficients}\end{equation} The functions $N^i_j(x,y)$,
defined on domains of induced local charts, are called the
\textit{local coefficients} of the nonlinear connection. The dual
basis of the adapted basis is $\{dx^i, \delta y^i=dy^i +
N^i_jdx^j\}$.

Every semispray $S$ determines a nonlinear connection. The
horizontal projector that corresponds to this nonlinear connection
is given by \cite{[grifone]}:
\begin{equation}
h(X)=\frac{1}{2}(Id-\mathcal{L}_SJ)(X)=\frac{1}{2}(X-[S,JX]+J[S,X]).
\label{hpojector}
\end{equation}
Local coefficients of the induced nonlinear connection are defined
on domains of induced local charts and they are given by
$N^i_j={\partial G^i}/{\partial y^j}$, \cite{[crampin]},
\cite{[grifone]}.

\section{Geometric structures on Lagrange space}

The presence of a regular Lagrangian on the tangent bundle $TM$
determines the existence of some geometric structures one can
associate to it such as: semispray, nonlinear connection and
symplectic structure.

Consider $L^{n}=(M, L)$ a Lagrange space. This means that $L:TM
\longrightarrow \mathbb{R}$ is differentiable of
$C^{\infty}$-class on $\widetilde{TM}$ and only continuous on the
null section. We also assume that $L$ is a regular Lagrangian. In
other words, the (0,2)-type, symmetric, d-tensor field with
components
\begin{equation} g_{ij}=\frac{1}{2}\frac{\partial^2 L}{\partial
y^i\partial y^j} \textrm{ has rank n on } \widetilde{TM}.
\label{lmetric}
\end{equation}

The Cartan 1-form $\theta_L$ of the Lagrange space can be defined
as follows:
\begin{equation}
\theta_L=J^*(d L)=d_JL=\frac{\partial L}{\partial y^i} dx^i.
\label{theta}
\end{equation}
For a vector field $X=X^i({\partial}/{\partial x^i}) +
Y^i({\partial}/{\partial y^i})$ on $TM$, the following formulae
are true:
\begin{equation}
\theta_L(X)=dL(JX)=d_JL(X)=(JX)(L)=\frac{\partial L}{\partial
y^i}X^i. \label{thetax}
\end{equation}

The Cartan 2-form $\omega_L$ of the Lagrange space can be defined
as follows:
\begin{equation}
\omega_L=d \theta_L = d(J^*(dL))=dd_JL= d\left(\frac{\partial
L}{\partial y^i}dx^i\right). \label{omega} \end{equation} In local
coordinates, the Cartan 2-form $\omega_L$ has the following
expression:
\begin{equation}
\omega_L =2g_{ij} dy^j \wedge d x^i +
\displaystyle\frac{1}{2}\left(\displaystyle\frac{\partial^2
L}{\partial y^i\partial x^j} - \displaystyle\frac{\partial^2
L}{\partial x^i \partial y^j}\right)d x^j\wedge d x^i.
\label{omegal}
\end{equation}
We can see from expression (\ref{omegal}) that the regularity of
the Lagrangian $L$ is equivalent with the fact that the Cartan
2-form $\omega_L$ has rank 2n on $\widetilde{TM}$ and hence it is
a symplectic structure on $\widetilde{TM}$.

The canonical semispray of the Lagrange space $L^n$ is the unique
vector field $S$ on $TM$ that satisfies the equation
\begin{equation}
i_S\omega_L = d(L-\mathbb{C}(L)). \label{isomega}
\end{equation}
The local coefficients $G^i$ of the canonical semispray $S$ are
given by the following formula:
\begin{equation}
G^i=\frac{1}{4}g^{ik}\left(\frac{\partial^2 L}{\partial y^k
\partial x^h} y^h-\frac{\partial L}{\partial x^k}\right).
\label{gi}
\end{equation}
Using the canonical semispray $S$ we can associate to a regular
Lagrangian $L$ a canonical nonlinear connection with the
horizontal projector given by expression (\ref{hpojector}) and the
local coefficients given by expression $N^i_j={\partial
G^i}/{\partial y^j}$.

The horizontal subbundle $NTM$ that corresponds to the canonical
nonlinear connection is a Lagrangian subbundle of the tangent
bundle $TTM$ with respect to the symplectic structure $\omega_L$.
This means that $\omega_L(hX,hY)=0$, $\forall X,Y, \in \chi(TM)$.
In local coordinates this implies the following expression for the
symplectic structure $\omega_L$:
\begin{equation}
\omega_L=2g_{ij}\delta y^j \wedge dx^i. \label{omegalocal}
\end{equation}

\section{Basic properties of a Lagrange space}

In this section we determine some basic properties for a
Lagrangian $L$ using the Cartan forms, canonical semispray and
nonlinear connection. If $h$ is the horizontal projector of the
canonical nonlinear connection of a regular Lagrangian we
determine a formula for the horizontal differential $d_hL$ of a
regular Lagrangian. From this formula we conclude that there are
non homogenous Lagrangians for which $d_hL\neq 0$.

\begin{proposition} \label{result1}
The following formulae regarding the Cartan 1-form $\theta_L$ of a
Lagrange space are true:
\begin{equation}
\begin{array}{l}
\iota_S\theta_L=\mathbb{C}(L), \vspace{2mm}\\
\mathcal{L}_S\theta_L = dL.
\end{array}
\label{ftheta}
\end{equation}
\end{proposition}
\begin{proof}
First formula (\ref{ftheta}) follows from the following
computation:
$$
\iota_S\theta_L=\theta_L(S)=\left(\frac{\partial L}{\partial
y^i}dx^i\right)\left(y^i\frac{\partial}{\partial x^i}-2G^i(x,y)
\frac{\partial}{\partial y^i}\right) = \frac{\partial L}{\partial
y^i} y^i =\mathbb{C}(L).
$$
If we differentiate first formula (\ref{ftheta}) we obtain
$d\iota_S\theta_L = d\mathbb{C}(L)$. Using the expression of the
Lie derivative $\mathcal{L}_S=d\iota_S + \iota_S d$ we obtain
$\mathcal{L}_S\theta_L -\iota_Sd\theta_L = d\mathbb{C}(L)$. Using
the defining formulae (\ref{omega}) for $\omega_L$ and
(\ref{isomega}) for the canonical semispray $S$ we obtain
$\mathcal{L}_S\theta_L=\iota_S\omega_L + d\mathbb{C}(L)=dL$ and
hence second formula (\ref{ftheta}) is true.
\end{proof}
\begin{theorem} \label{result2}
Consider $h$ the horizontal projector (\ref{hpojector}) of a
Lagrange space. We have the following formula for the horizontal
differential operator $d_h$:
\begin{equation}
d_hL=\frac{1}{2}d_J(S(L)). \label{dhl}
\end{equation}
In local coordinates, formula (\ref{dhl}) is equivalent with the
following expression for the horizontal covariant derivative of
the Lagrangian $L$:
\begin{equation}
L_{|i}:=\frac{\delta L}{\delta x^i} = \frac{1}{2}
\frac{\partial}{\partial y^i}(S(L)). \label{dhllocal}
\end{equation}
\end{theorem}
\begin{proof}
In order to prove (\ref{dhl}) we have to show that for every $X\in
\chi(TM)$, we have that
$$ (d_hL)(X):=dL(hX)= \frac{1}{2}(JX)(S(L)).$$
Using second formula (\ref{ftheta}) we obtain
$$
\begin{array}{ll}
0 & =(\mathcal{L}_S\theta_L - dL)(X) =
S\theta_L(X)-\theta_L[S,X]-dL(X) \vspace{2mm}\\
& = S((JX)(L))-J[S,X](L)-dL(X)  \vspace{2mm}\\
& = [S,JX](L) + (JX)(S(L)) - J[S,X](L) - dL(X) \vspace{2mm}\\
& = (JX)(S(L)) - dL(X - [S,JX] + J[S,X]) \vspace{2mm}\\
& = (JX)(S(L))-dL(2hX).
\end{array}
$$
Consequently, formula (\ref{dhl}) is true.

Due to the linearity of the operators involved in formula
(\ref{dhl}) we have that formulae (\ref{dhl}) and (\ref{dhllocal})
are equivalent. However, we give here an independent proof that
formula (\ref{dhllocal}) is true. The right hand side of formula
(\ref{dhllocal}) can be expressed as follows:
\begin{equation}
\begin{array}{rl}
\displaystyle\frac{\partial}{\partial y^i}(S(L)) =&
\displaystyle\frac{\partial}{\partial
y^i}\left(\displaystyle\frac{\partial L}{\partial x^j}y^j -
2\displaystyle\frac{\partial L}{\partial y^j}G^j\right) \vspace{2mm}\\
= & \displaystyle\frac{\partial L}{\partial x^i} +
\displaystyle\frac{\partial^2 L}{\partial y^i \partial x^j}y^j -
4g_{ij} G^j - 2 \displaystyle\frac{\partial L}{\partial y^j}N^j_i
\vspace{2mm} \\
= & 2\left(\displaystyle\frac{\partial L}{\partial x^i} -
\displaystyle\frac{\partial L}{\partial y^j}N^j_i\right) =
2\displaystyle\frac{\delta L}{\delta x^i}=2L_{|i}.
\end{array}
\end{equation}
In the above calculation we did use expression (\ref{gi}) for the
local coefficients $G^i$ of the canonical semispray $S$.
\end{proof}
\begin{corollary} \label{corollary1}
If the Lagrangian $L$ is constant along the solution curves of the
Euler-Lagrange equations then $L$ is constant along the horizontal
curves of the canonical nonlinear connection.
\end{corollary}
\begin{proof}
It is immediate from (\ref{dhl}) that $S(L)=0$ implies $d_hL=0$.
\end{proof}

\begin{proposition} \label{result3}
Consider a regular Lagrangian $L$ that is homogeneous of order k,
$k\neq 1$, with respect to $y$. The horizontal differential of the
Lagrangian $L$ vanishes, which means that $d_hL=0.$
\end{proposition}
\begin{proof}
Using Euler theorem for homogeneous functions we have that the
Lagrangian $L$ is homogeneous of order k if and only if
$\mathbb{C}(L)=({\partial L}/{\partial y^i})y^i=kL.$ Hence,
formula (\ref{isomega}), which defines the canonical semispray of
the Lagrangian $L$, can be written as $\iota_S\omega_L=(1-k)dL$.
This implies that $(1-k)S(L)=(1-k)dL(S)=\omega_L(S,S)=0$. Using
expression (\ref{dhl}) we obtain $d_hL=0$.
\end{proof}

Regular Lagrangians that are second order homogeneous with respect
to $y^i$ are encountered in Finsler geometry. For a Finsler space
its geodesics with the arclength parameterization coincide with
the integral curves of the canonical semispray, which are
horizontal curves with respect to the nonlinear connection. The
property $d_hL=0$ tells us that the Lagrangian $L$ is constant
along the horizontal curves of the nonlinear connection and hence
it is constant along the geodesics of the space.

\begin{theorem} \label{result4}
There exist regular Lagrangians $L$ for which $d_hL\neq 0$.
\end{theorem}
\begin{proof}
Let $a_{ij}(x)$ be a Riemannian structure on the base manifold $M$
and $\varphi$ a function on the manifold $M$ that will be
determined latter. The function $L'(x,y)=a_{ij}(x)y^iy^j$ is a
second order homogeneous regular Lagrangian on $TM$. The following
function is also a regular Lagrangian:
\begin{equation}
L(x,y)=L'(x,y)+\frac{\partial \varphi}{\partial x^i}(x)y^i.
\label{ll'}
\end{equation}
It is a straightforward calculation to check that $L'$ and
$L=L'+\varphi^c$ have the same metric tensor $a_{ij}$ and
consequently since $L'$ is a regular Lagrangian so is $L$. One can
also check that the Cartan 2-forms of the two Lagrangians
coincide, which means that $\omega_{L'}=\omega_L$. Therefore using
expression (\ref{isomega}) the two Lagrangians have the same
semispray $S$ and using (\ref{hpojector}) they have the same
nonlinear connection with the same horizontal projector $h$. Local
coefficients of the canonical semispray $S$ are
$2G^i(x,y)=\gamma^i_{jk}(x)y^jy^k$, where $\gamma^i_{jk}(x)$ are
the Christoffel's symbols of the second kind for the Riemannian
metric $a_{ij}$. According to formula (\ref{dhl}) and Proposition
\ref{result3} we have
\begin{equation}
\begin{array}{ll}
(d_hL)(X) & =\displaystyle\frac{1}{2}(JX)(S(L'+\varphi^c)) =
\displaystyle\frac{1}{2}(JX)(S(\varphi^c)) \vspace{2mm}\\
&= \displaystyle\frac{1}{2}\left(2\displaystyle\frac{\partial^2
\varphi}{\partial x^i\partial x^j}y^j - 2
\displaystyle\frac{\partial \varphi}{\partial x^j}
\displaystyle\frac{\partial G^j}{\partial y^i}\right)X^i
\vspace{2mm}\\
& = \left(\displaystyle\frac{\partial^2 \varphi}{\partial
x^i\partial x^j} - \displaystyle\frac{\partial \varphi}{\partial
x^k} \gamma^k_{ij}\right)y^jX^i,
\end{array}
\label{contra}
\end{equation}
where $X=X^i({\partial}/{\partial x^i}) + Y^i({\partial}/{\partial
y^i})$ is a vector field on $TM$. At this moment we can choose a
function $\varphi$ on $M$ such that
\begin{equation}
\frac{\partial^2 \varphi}{\partial x^i\partial x^j} \neq
\frac{\partial \varphi}{\partial x^k} \gamma^k_{ij}.\label{neqphi}
\end{equation}
For example we can choose either the function $\varphi$ to be
linear in $x$ and $a_{ij}$ a non flat Riemannian metric or we can
choose the Riemannian metric $a_{ij}$ to be flat and $\varphi$ a
nonlinear function. For the first example, the left hand side of
formula (\ref{neqphi}) is zero, while the right hand side is not.
For the second example we have that the right hand side of formula
(\ref{neqphi}) is zero, while the left hand side is not zero.

With a function $\varphi$ and a Riemannian metric $a_{ij}$ that
satisfy (\ref{neqphi}), we have that for the Lagrangian $L$
defined by (\ref{ll'}) the horizontal differential is not zero,
which means that $d_hL\neq 0$.
\end{proof}

We want to mention that for the Lagrangian (\ref{ll'}) we have
that $S(L)=0$ if and only if $d_h(L)=0$. However we don't know if
this result is true for an arbitrary Lagrangian $L$. Hence we
don't know if the converse of Corollary \ref{corollary1} is true
which means that $d_hL=0$ implies $S(L)=0$?

\vspace*{3mm}

\noindent\textbf{Acknowledgment} This work has been partially
supported by CNCSIS grant from the Romanian Ministry of Education.


\begin{thebibliography}{n}

\bibitem[Cra71]{[crampin]} Crampin, M.: On horizontal distributions
on the tangent bundle of a differentiable manifold. J. London
Math. Soc. \textbf{2} (3), 178--182 (1971).

\bibitem[Gri72]{[grifone]}
Grifone, J.: Structure presque-tangente et connexions I. Ann.
Inst. Henri Poincare. \textbf{22} (1), 287--334 (1972).

\bibitem[Kle82]{[klein]}
Klein, J.: Geometry of Sprays. Lagrangian case. Principle of Least
curvature. IUTAM-ISIMM Symposium on Analytical Mechanics, Torino,
177--196 (1982).

\bibitem[Sze93]{[szenthe1]}
Szenthe, J.: On a basic property of Lagrangians. Publ. Math.
Debrecen. \textbf{42} (3--4), 247--251 (1993).

\bibitem[Sze96]{[szenthe2]}
Szenthe, J.: Symmetries of Lagrangian fields are homoteties in the
homogeneous case. Proceedings of the Conference on Differential
Geometry and its Applications, Brno, 321--328 (1996).

\bibitem[SM93]{[szilasi]}
Szilasi, J. Muzsnay, Z.: Nonlinear connections and the problem of
metrizability. Publ. Math. Debrecen. \textbf{42} (1--2), 175--192
(1993).


\end{thebibliography}
\end{document}